\documentclass[11pt,showlabels]{article1}
\usepackage{amssymb}
\usepackage{epsfig}
\usepackage{amsmath}
\title{\bf Optimal Control and
Stabilization Problem for Discrete-time Markov Jump Systems
with Indefinite Weight Costs
 \thanks{
 This work is supported by the National Natural Science
Foundation of China under Grants  61573221,
61633014, 61473134. 
$^{*}$Corresponding author: Huanshui Zhang. Email: hszhang@sdu.edu.cn}
}

\author{Hongdan Li$^a$,\ Chunyan Han$^b$,\ Huanshui Zhang$^{a,*}$
\ \\
\\
\ \ \ $^a$ School of Control Science and Engineering, Shandong
University, \\Jinan Shandong 250061, China
\\
$^b$  School of Electrical Engineering, University of Jinan, \\Jinan Shandong 250022, China }

\begin{document}
\baselineskip 16pt
\date{}
  \maketitle
\begin{abstract} It is well known that stability is the most fundamental nature with regard to a control system, in view of this, the stabilization becomes an inevitable control problem. This article mainly discusses the optimal control and
 stabilization problem for discrete-time systems involving Markov jump and multiplicative noise. The state and control weighting matrices in the cost function are allowed to be indefinite.
By solving the forward-backward stochastic difference equations with Markov jump (FBSDEs-MJ) derived from the maximum principle, we conclude that the necessary and sufficient conditions of the solvability of indefinite optimal control problem in finite-horizon, whose method is different from most previous works \cite{13}, etc. Furthermore,
 necessary and sufficient conditions that stabilize the Markov jump discrete-time systems in the mean square sense with indefinite weighting matrices in the cost are first developed under the basic assumption of exactly observable, which is different from the previous works \cite{12}, \cite{14} where an additional assumption of stabilization of systems is made.

The key points of this article can be summed up as that an analytic solution to FBSDEs-MJ which makes the optimal controller to be explicitly expressed  and the method of transformation, i.e., the stabilization problem of indefinite case is boiled down to a definite one whose stabilization is expressed by defining Lyapunov function via the optimal cost subject to a new algebraic Riccati equation involving Markov jump (NGARE-MJ).

\bigskip

\end{abstract}

\pagestyle{plain} \setcounter{page}{1}
\section{Introduction}

 There are many factors to give rise to abrupt changes such as abrupt environmental disturbances, component failures
or repairs and these changes often occur in many control systems, for instance, economic systems and aircraft control systems. This phenomenon can be modeled as Markov jump linear systems (MJLS). Owing to its widely application in practice, in recent years,  the subject of MJLS is by now huge and is growing rapidly, see \cite{1}-\cite{7}, and reference therein. Seeing that the importance of the linear quadratic (LQ) control and stabilization problem in the study of control system, there are also many results about these problems with Markov jump. \cite{8} considered the optimal control problems for  discrete-time linear systems subject to Markov jump with two cases that the one without noise and the other with an additive noise in model. In \cite{9}, they illustrated the equivalence between the stability of the optimal control and positiveness of  the coupled algebraic Riccati equation via the concept of weak detectability.

It is noteworthy that all the above results are obtained under the common assumption that the weighting matrices of state and control in the quadratic performance index are required to be positive semi-definite even positive definite. However, when the weighting matrices have the requirement of symmetry only, the stochastic LQ problem may be still well posed. This case is called indefinite stochastic problem which often appear in economic fields such as portfolio selection problem.
 As regard to the problem, \cite{10} first considered the well-posedness of indefinite LQ problem for continuous-time system. Since then,  \cite{11} and \cite{12} investigated an indefinite stochastic LQ control problem for  continuous-time linear systems subject to Markov jump in finite and infinite time horizon, respectively. As to discrete-time linear systems with Markov jump, \cite{13} derived the necessary and sufficient condition for the well posedness of the indefinite LQ problem and the optimal control law were given in terms of a set
of coupled generalized Riccati difference equations interconnected with a set of coupled linear recursive equations. In \cite{14}, under the assumption that the system is mean square stabilizable, they gave the sufficient conditions for the existence of the maximal solution, and necessary and sufficient conditions for the existence of the mean square stabilizing solution for the generalized coupled algebraic Riccati equations.

What is worth considering is that there were few articles to consider the stabilization problem for discrete-time linear systems subject to Markov jump with indefinite weighting matrices in cost function. However, for most situation, stability is the precondition for the normal operation of a control system. More concretely,  the study of infinite optimal control problem is of significance only under the condition that the system is stabilizable, see reference \cite{14}-\cite{17}. Actually, there have been some results of stabilization about definite problem,  such as, \cite{18} and \cite{19} considered the stabilization problem with positive semi-definite and positive definite weighting matrices for stochastic systems involving multiplicative noises and input delay in the case of discrete-time and continuous-time, respectively. \cite{20} provided the stabilization for discrete-time mean-field systems with positive semi-definite weighting matrices. However, the above conclusions are either obtained under the definite condition, or the stabilization problem is not given in indefinite case.  In view of this, in this paper, we mainly consider the optimal control and
 stabilization problem for discrete-time systems subject to Markov jump and multiplicative noise, where the state and control weighting matrices in the cost function are allowed to be indefinite.

The main contribution of this paper can be summed up as that an optimal controller is explicitly shown by a generalized difference Riccati equation with Markov jump (GDRE-MJ) which is derived from the solution to the FBSDEs-MJ, which is a new method compared with the previous works studied the linear quadratic optimal problem involving Markov jump. Secondly, under the preconditions that a set involving linear matrix inequality  and kernel restriction is empty and system is exactly observable, we derive the existence of the maximum solution to GARE-MJ by discussing the convergence of  the associated GDRE-MJ.  This precondition is more easily verified compared with the requirement in \cite{14} about some operators and even spectral radius. The main results in \cite{14} are obtained based on the stabilization of the system and there is no discussion about the stabilization. Therefore,  another contribution in our paper is that the conclusion of stabilization for discrete-time systems subject to Markov jump and multiplicative noise with indefinite weighting matrices is expressed for the first time under the basic assumption of exactly observable, which is different from the previous works.

 The rest of this article is made up of the following sections. Section 2 mainly provides some results about optimal control with finite horizon. The conclusion of stabilization will be shown in section 3. We will give a numerical example in section 4 to further illustrate the correctness of the conclusion. And in section 5 we will make a summary.

The related notations in this article are expressed as follows:  \\
${\mathbb{R}}^n$ : the $n$-dimensional Euclidean space;\\
$\mathbb{R}^{m\times n}$ : the norm bounded linear space of all
$m\times n$ matrices;\\
$Y'$ : the transposition of $Y$;\\
$Y\geq 0 (Y>0)$: the symmetric matrix $Y\in \mathbb{R}^{n\times n}$ is positive semi-definite(positive definite);\\
$Y^{\dagger}$ : the Moore-Penrose pseudo-inverse of $Y$;\\
$\mathbf{Ker}(Y)$ : the kernel of a matrix $Y$;\\
$(\Omega,\mathcal{G}, \mathcal{G}_{k}, \mathcal{P})$: a complete probability space with the $\sigma$-field generated by $\{x(0),\theta(0),\cdots, x(k),\theta(k)\}$;\\
$E[\cdot|\mathcal{G}_{k}]$:  the conditional expectation with respect
to $\mathcal{G}_{k}$ and $\mathcal{G}_{-1}$ is understood as $\{\emptyset,\Omega\}$.

\section{Preliminaries }


\setcounter{equation}{0}
Considering  the following discrete-time  Markov jump linear system with multiplicative noise:
\begin{eqnarray}
x(k+1)&=&(A_{\theta(k)}+B_{\theta(k)}\omega(k))x(k)+(C_{\theta(k)}+D_{\theta(k)}\omega(k))u(k),\label{x1}
\end{eqnarray}
where $x(k)\in {\mathbb{R}}^n$ denotes the state, $u(k)\in {\mathbb{R}}^m$  denotes  control process and $\omega(k)$ is scalar valued random white
noise with zero mean and variance $\sigma^{2}$. $\theta(k)$ is a discrete-time Markov chain with finite state space $\{1,2,\cdots,L\}$ and transition probability $\rho_{i,j}=\mbox{P}(\theta(k+1)=j|\theta(k)=i)(i,j=1,2,\cdots,L)$. We set $\pi_{i}(k)=\mbox{P}(\theta(k)=i)(i=1,2,\cdots,L)$, while $A_{i}, B_{i}, C_{i},D_{i}({i}=1,\cdots,L)$ are matrices of appropriate dimensions. The initial value $x_0$ is known. We assume that $\theta(k)$ is independent of $x_0$.

The quadratic cost subject to  system (\ref{x1}) with infinite horizon is given by
\begin{eqnarray}
J&=&\mbox{E}\bigg\{\sum_{k=0}^{\infty}[x(k)'Q_{\theta(k)}x(k)+u(k)'R_{\theta(k)}u(k)]\bigg\},\label{J2}
\end{eqnarray}
where $Q_{\theta(k)},R_{\theta(k)}$ are just symmetric matrices.

The following problem will be mainly discussed in this paper, i.e.,

{\bf Problem 1} \ \ Find the $\mathcal{G}_{k}$-measurable controller $u(k)=F_{\theta(k)}x(k)$ with constant
matrix gain $F_{\theta(k)}$ to stabilize (\ref{x1}) while minimizing (\ref{J2}).\\

 While for the convenience of discussing the above problem, we will first introduce some associated results about the cost function with finite horizon as the following description.
\begin{eqnarray}
J_N&=&\mbox{E}\bigg\{\sum_{k=0}^N\left[x(k)'Q_{\theta(k)}x(k)+u(k)'R_{\theta(k)}u(k)\right]\nonumber\\
&&+x(N+1)'P_{\theta(N+1)}x(N+1)\bigg\},\label{J3}
\end{eqnarray}
where $N>0$ is an integer, $x(N+1)$ is the terminal state, $P_{\theta(N+1)}$ reflects the penalty on the terminal state, the matrix functions $R_{\theta(k)}$ and $Q_{\theta(k)}$ are symmetric matrices.

As to the case of finite horizon, we will discuss the Problem$^{\ast}$, i.e.,

{\bf Problem$^{\ast}$} \ \ Find a $\mathcal{G}_{k}$-measurable controller $u(k)$ to minimize (\ref{J3}) subject to (\ref{x1}).

On the ground of the indefiniteness of weighting matrices, the above problem may be ill-posed. Hence, we should introduce next definitions and lemmas.

{\bf Definition 1}: Problem$^{\ast}$ is called well posed if
$\inf\limits_{u_{0},\cdots, u_{N}}J_N>-\infty$
for any random variables $x_{0}$.

{\bf Definition 2}: Problem$^{\ast}$ is called solvable if there exists an admissible control $(u^{\ast}_{0},\cdots, u^{\ast}_{N})$ such that (\ref{J3}) is minimized for any $x_{0}$.

{\bf Remark 1}: From Theorem 4.3 in \cite{21}, the equivalence between the well-posedness  and the solvability of Problem$^{\ast}$ can be obtained.

The following lemmas are about some properties of the pseudo inverse matrix.

{\bf Lemma 1 \cite{22} } \ \ Let a symmetric matrix $S$ be given. Then
\begin{itemize}
              \item [(i)] $S^{\dagger}=S^{\dagger'}$;
              \item [(ii)]$S\geq0$ if and only if $S^{\dagger}\geq0$;
              \item [(iii)]$SS^{\dagger}=S^{\dagger}S$.
            \end{itemize}

{\bf Lemma 2 \cite{22} } \  Let a matrix $M\in \mathbb{R}^{m\times n}$ be given. Then there exists a unique matrix $M^{\dagger}\in \mathbb{R}^{n\times m}$ such that
\begin{itemize}
              \item [(i)] $MM^{\dagger}M=M,M^{\dagger}MM^{\dagger}=M^{\dagger}$;
              \item [(ii)]$(MM^{\dagger})'=MM^{\dagger},(M^{\dagger}M)'=M^{\dagger}M$.
            \end{itemize}

{\bf Lemma 3 (Extended Schur’s Lemma \cite{22})} \ \ Let $M=M'$, $N$, and $R=R'$ be given matrices with appropriate dimensions. Then the following conditions are equivalent:
\begin{itemize}
\item [(i)]  $M-NR^{\dagger}N'\geq 0, R\geq 0, N(I-RR^{\dagger})=0$;
\item [(ii)]   $\left[
  \begin{array}{cc}
    M & N\\
    N' & R\\  
  \end{array}
\right]\geq 0 $;
\item [(iii)]   $\left[
  \begin{array}{cc}
    R & N'\\
    N & M\\  
  \end{array}
\right]\geq 0$.
\end{itemize}

Due to the dependence of $\theta(k)$ on its past values, an extended version of the stochastic maximum principle which is suitable for the MJLS (1) is established in the sequel.

{\bf Lemma 4 (Maximum Principle involving Markov Jump)}  \ \
According to the linear system (\ref{x1}) and the performance index (\ref{J3}). If the linear quadratic  problem $\min J_N$ is solvable, then the optimal ${\cal{G}}_{k}$-measurable control $u(k)$ satisfies the following equilibrium condition
\begin{eqnarray}
0=\mbox{E}[(C_{\theta(k)}+D_{\theta(k)}\omega(k))'\lambda_k+R_{\theta(k)}u(k)|{\cal{G}}_{k}],k=0,\cdots,N,\label{f4}
\end{eqnarray}
where the costate $\lambda_k$ satisfies the following equation
\begin{eqnarray}
\lambda_N&=&\mbox{E}[P_{\theta(N+1)}x(N+1)|{\cal{G}}_N],\label{f5}\\
\lambda_{k-1}&=&\mbox{E}[(A_{\theta(k)}+B_{\theta(k)}\omega(k))'\lambda_k+Q_{\theta(k)}x(k)|{\cal{G}}_{k-1}],k=0,\cdots,N,\label{ff6}
\end{eqnarray}

together the costate equation (\ref{f5})-(\ref{ff6}) with state equation (\ref{x1}), the FBSDEs-MJ is established, which play a vital role in this paper.

\emph{Proof}. Similar to the derivation for Maximum Principle (MP) as in \cite{18},\cite{24}, the MP (\ref{f4})-(\ref{ff6}) follows directly, the aforementioned conclusion can be derived using an analogous step, so its proof is omitted.

Now we will show the following theorem which is expressed the result of Problem$^{\ast}$.

{\bf Theorem 1}
Problem$^{\ast}$ is solvable if and only if the following  generalized difference Riccati equations with Markov jump
\begin{eqnarray}
\left\{
\begin{array}{lll}
P_{i}(k)=A_{i}'(\sum_{j=1}^L\rho_{i,j}P_{j}(k+1))A_{i}
+\sigma^{2}B_{i}'(\sum_{{j}=1}^L\rho_{i,j}P_{j}(k+1))B_{i}+Q_{i}\\
-M_{i}(k)'\Upsilon_{i}(k)^{\dagger}M_{i}(k),\\
\Upsilon_{i}(k)\Upsilon_{i}(k)^{\dagger}M_{i}(k)-M_{i}(k)=0,\\
\Upsilon_{i}(k)\geq0,\label{f7}
\end{array}
\right.
\end{eqnarray}
in which
\begin{eqnarray}
\Upsilon_{i}(k)&=&C_{i}'(\sum_{{j}=1}^L\rho_{i,j}P_{j}(k+1))C_{i}
+\sigma^{2}D_{i}'(\sum_{{j}=1}^L\rho_{i,j}P_{j}(k+1))D_{i}+R_{i},\label{f8}\\
M_{i}(k)&=&C_{i}'(\sum_{j=1}^L\rho_{i,j}P_{j}(k+1))A_{i}
+\sigma^{2}D_{i}'(\sum_{{j}=1}^L\rho_{i,j}P_{j}(k+1))B_{i},\label{f9}
\end{eqnarray}
has a solution. If this condition is satisfied, the analytical solution to the optimal control can be given as
\begin{eqnarray}
u^{\ast}(k)=-\Upsilon_{i}(k)^{\dagger}M_{i}(k){x}(k), i=1,\cdots,L,\label{f10}
\end{eqnarray}
for $k=N,\cdots,0$.
The corresponding optimal performance index is given by
\begin{eqnarray}
J^{\ast}_N=\mbox{E}[x(0)'P_{\theta(0)}(0)x(0)].\label{f11}
\end{eqnarray}
The relationship of the costate $\lambda_{k-1}$ and the state $x(k)$ is given as
\begin{eqnarray}
\lambda_{k-1}=(\sum_{j=1}^L\rho_{i,j}P_{j}(k))x(k), i=1,\cdots,L.\label{f12}
\end{eqnarray}

\emph{Proof}. (Necessity) Assume that Problem$^{\ast}$ is solvable, we will investigate that there exist symmetric matrices $P_{i}(0), \cdots, P_{i}(N)$, $i=1,\cdots,L$ satisfying the GDRE-MJ (\ref{f7}) by induction. To this end, we first set the following formula as
\begin{eqnarray}
\underline{J}(k)&=&\inf\limits_{u_{k},\cdots, u_{N}}\mbox{E}\bigg[\sum_{{i}=k}^N(x(i)'Q_{\theta(i)}x(i)+u(i)'R_{\theta(i)}u(i))\nonumber\\
&&+x(N+1)'P_{\theta(N+1)}x(N+1)|{\cal{G}}_{k-1}\bigg].\label{f13}
\end{eqnarray}
It is obvious to know that for any $k_{1}<k_{2}$, when $\underline{J}(k_{1})$ is finite then $\underline{J}(k_{2})$ is also finite by the stochastic optimality principle. Since Problem$^{\ast}$ is supposed to be solvable, we can see that $\underline{J}(k)$ is finite for any $0\leq k\leq N$.

Firstly, we let $k=N$, from system (\ref{x1}), we know that
\begin{eqnarray*}
\underline{J}(N)&=&\inf\limits_{ u_{N}}\mbox{E}\Bigg\{x(N)'Q_{\theta(N)}x(N)+u(N)'R_{\theta(N)}u(N)\nonumber\\
&&+[(A_{\theta(N)}+B_{\theta(N)}\omega(N))x(N)+(C_{\theta(N)}+D_{\theta(N)}\omega(N))u(N)]'P_{\theta(N+1)}\\
&&\cdot[(A_{\theta(N)}+B_{\theta(N)}\omega(N))x(N)+(C_{\theta(N)}+D_{\theta(N)}\omega(N))u(N)]|{\cal{G}}_{N-1}\Bigg\}\\
&=&\inf\limits_{ u_{N}}\mbox{E}\Bigg\{x(N)'[Q_{\theta(N)}+A_{\theta(N)}'P_{\theta(N+1)}A_{\theta(N)}
+\sigma^{2}B_{\theta(N)}'P_{\theta(N+1)}B_{\theta(N)}]x(N)\\
&&+2x(N)'[A_{\theta(N)}'P_{\theta(N+1)}C_{\theta(N)}
+\sigma^{2}B_{\theta(N)}'P_{\theta(N+1)}D_{\theta(N)}]u(N)\\
&&+u(N)'[R_{\theta(N)}+C_{\theta(N)}'P_{\theta(N+1)}C_{\theta(N)}
+\sigma^{2}D_{\theta(N)}'P_{\theta(N+1)}D_{\theta(N)}]u(N)|{\cal{G}}_{N-1}\Bigg\}\\
&=&\inf\limits_{ u_{N}}\mbox{E}\Bigg\{x(N)'[Q_{i}+A_{i}'(\sum_{{j}=1}^L\rho_{i,j}P_{j}(N+1))A_{i}
+\sigma^{2}B_{i}'(\sum_{{j}=1}^L\rho_{i,j}P_{j}(N+1))B_{i}]x(N)\\
&&+2x(N)'[A_{i}'(\sum_{{j}=1}^L\rho_{i,j}P_{j}(N+1))C_{i}
+\sigma^{2}B_{i}'(\sum_{{j}=1}^L\rho_{i,j}P_{j}(N+1))D_{i}]u(N)\\
&&+u(N)'[R_{i}+C_{i}'(\sum_{{j}=1}^L\rho_{i,j}P_{j}(N+1))C_{i}
+\sigma^{2}D_{i}'(\sum_{{j}=1}^L\rho_{i,j}P_{j}(N+1))D_{i}]u(N)|{\cal{G}}_{N-1}\Bigg\}
\end{eqnarray*}
By Lemma 4.3 in \cite{21} and the finiteness of $\underline{J}(N)$, it yields that there indeed exist symmetric matrix $P_{i}(N)$ satisfying
\begin{eqnarray*}
\underline{J}(N)=\mbox{E}[x(N)'P_{i}(N)x(N)],
\end{eqnarray*}
and furthermore,
\begin{eqnarray}
&&P_{i}(N)=A_{i}'(\sum_{j=1}^L\rho_{i,j}P_{j}(N+1))A_{i}
+\sigma^{2}B_{i}'(\sum_{{j}=1}^L\rho_{i,j}P_{j}(N+1))B_{i}+Q_{i}\nonumber\\
&&\hspace{10mm}-M_{i}(N)'\Upsilon_{i}(N)^{\dagger}M_{i}(N),\label{f14}\\
&&\Upsilon_{i}(N)\Upsilon_{i}(N)^{\dagger}M_{i}(N)-M_{i}(N)=0,\label{f15}\\
&&\Upsilon_{i}(N)\geq0,\label{f16}
\end{eqnarray}
in which
\begin{eqnarray}
\Upsilon_{i}(N)&=&C_{i}'(\sum_{{j}=1}^L\rho_{i,j}P_{j}(N+1))C_{i}
+\sigma^{2}D_{i}'(\sum_{{j}=1}^L\rho_{i,j}P_{j}(N+1))D_{i}+R_{i},\label{f17}\\
M_{i}(N)&=&C_{i}'(\sum_{j=1}^L\rho_{i,j}P_{j}(N+1))A_{i}
+\sigma^{2}D_{i}'(\sum_{{j}=1}^L\rho_{i,j}P_{j}(N+1))B_{i}.\label{f18}
\end{eqnarray}

The optimal controller $u(N)$ will be calculated from (\ref{x1}), (\ref{f4}) and (\ref{f5}).
\begin{eqnarray}
0&=&E[(C_{\theta(N)}+D_{\theta(N)}\omega(N))'\lambda(N)+R_{\theta(N)}u(N)|\mathcal{G}_{N}]\nonumber\\
&=&E[(C_{\theta(N)}+D_{\theta(N)}\omega(N))'E[P_{\theta(N+1)}x(N+1)|{\cal{G}}_N]+R_{\theta(N)}u(N)|\mathcal{G}_{N}]\nonumber\\
&=&E[(C_{\theta(N)}+D_{\theta(N)}\omega(N))'\sum_{j=1}^{L}\rho_{i,j}P_{j}(N+1)x(N+1)|{\cal{G}}_N]+R_{\theta(N)}u(N)|\mathcal{G}_{N}]\nonumber\\
&=&\bigg[C_{i}'(\sum_{j=1}^L\rho_{i,j}P_{j}(N+1))A_{i}
+\sigma^{2}D_{i}'(\sum_{{j}=1}^L\rho_{i,j}P_{j}(N+1))B_{i}\bigg]x(N)\nonumber\\
&&+\bigg[C_{i}'(\sum_{{j}=1}^L\rho_{i,j}P_{j}(N+1))C_{i}
+\sigma^{2}D_{i}'(\sum_{{j}=1}^L\rho_{i,j}P_{j}(N+1))D_{i}+R_{i}\bigg]u(N).\label{f19}
\end{eqnarray}

So, from (\ref{f17}) and (\ref{f18}), we have that
\begin{eqnarray}
u(N)=-\Upsilon_{i}(N)^{\dagger}M_{i}(N)x(N),\label{f20}
\end{eqnarray}
which is as (\ref{f10}) in the case of $k=N$.

As to $\lambda_{N-1}$, from (\ref{x1}), (\ref{f5}), (\ref{ff6}) and (\ref{f20}), it yields that
\begin{eqnarray}
\lambda_{N-1}&=&\mbox{E}[(A_{\theta(N)}+B_{\theta(N)}\omega(N))'\lambda_N+Q_{\theta(N)}x(N)|{\cal{G}}_{N-1}]\nonumber\\
&=&\mbox{E}[(A_{\theta(N)}+B_{\theta(N)}\omega(N))'E[P_{\theta(N+1)}x(N+1)|{\cal{G}}_N]+Q_{\theta(N)}x(N)|{\cal{G}}_{N-1}]\nonumber\\
&=&\mbox{E}\bigg[A_{i}'(\sum_{j=1}^L\rho_{i,j}P_{j}(N+1))A_{i}+B_{i}'(\sum_{j=1}^L\rho_{i,j}P_{j}(N+1))B_{i}+Q_{i}\nonumber\\
&&-M_{i}(N)'\Upsilon_{i}(N)^{\dagger}M_{i}(N)|{\cal{G}}_{N-1}\bigg]{x}(N)\nonumber\\
&=&(\sum_{i=1}^L\rho_{s,i}P_{i}(N))x(N),s=1,\cdots,L,
\end{eqnarray}
which is satisfied (\ref{f12}) with $k=N$.

Now we assume that GDRE-MJ (\ref{f7}) has a solution $P_{i}(m)$, $k+1\leq m \leq N$ and satisfying $\underline{J}(m)=\mbox{E}[x(m)'P_{i}(m)x(m)]$ and $u(m)$, $\lambda(m-1)$ are as (\ref{f10}), (\ref{f12}), respectively, thus for $k$, we have
\begin{eqnarray*}
\underline{J}(k)&=&\inf\limits_{ u_{k}}\mbox{E}\left[x(k)'Q_{\theta(k)}x(k)+u(k)'R_{\theta(k)}u(k)+\underline{J}(k+1)|{\cal{G}}_{k-1}\right]\\
&=&\inf\limits_{ u_{k}}\mbox{E}\left\{x(k)'Q_{\theta(k)}x(k)+u(k)'R_{\theta(k)}u(k)+\mbox{E}[x(k+1)'P_{i}(k+1)x(k+1)]|{\cal{G}}_{k-1}\right\}\\
&=&\inf\limits_{ u_{k}}\mbox{E}\bigg\{x(k)'[Q_{i}+A_{i}'(\sum_{{j}=1}^L\rho_{i,j}P_{j}(k+1))A_{i}
+\sigma^{2}B_{i}'(\sum_{{j}=1}^L\rho_{i,j}P_{j}(k+1))B_{i}]x(k)\\
&&+2x(k)'[A_{i}'(\sum_{{j}=1}^L\rho_{i,j}P_{j}(k+1))C_{i}
+\sigma^{2}B_{i}'(\sum_{{j}=1}^L\rho_{i,j}P_{j}(k+1))D_{i}]u(k)\\
&&+u(k)'[R_{i}+C_{i}'(\sum_{{j}=1}^L\rho_{i,j}P_{j}(k+1))C_{i}
+\sigma^{2}D_{i}'(\sum_{{j}=1}^L\rho_{i,j}P_{j}(k+1))D_{i}]u(k)|{\cal{G}}_{k-1}\bigg\}.
\end{eqnarray*}
Similarly, from Lemma 4.3 in \cite{21} and the finiteness of $\underline{J}(k)$, we can obtain that there exist $P_{i}(k)$ satisfying GDRE-MJ (\ref{f7}). Furthermore, $\underline{J}(k)=\mbox{E}[x(k)'P_{i}(k)x(k)]$. From now on by mathematical induction we obtain that GDRE-MJ (\ref{f7}) exists a solution.

In the case that GDRE-MJ (\ref{f7}) exists a solution and the inductive hypothesis, the optimal controller $u(k)$ can be obtained from (\ref{x1}) and (\ref{f4}).
\begin{eqnarray}
0&=&E[(C_{\theta(k)}+D_{\theta(k)}\omega(k))'\lambda(k)+R_{\theta(k)}u(k)|\mathcal{G}_{k}]\nonumber\\
&=&E[(C_{\theta(k)}+D_{\theta(k)}\omega(k))'(\sum_{j=1}^L\rho_{i,j}P_{j}(k+1))x(k+1)+R_{\theta(k)}u(k)|\mathcal{G}_{k}]\nonumber\\
&=&\bigg[C_{i}'(\sum_{j=1}^L\rho_{i,j}P_{j}(k+1))A_{i}
+\sigma^{2}D_{i}'(\sum_{{j}=1}^L\rho_{i,j}P_{j}(k+1))B_{i}\bigg]x(k)\nonumber\\
&&+\bigg[C_{i}'(\sum_{{j}=1}^L\rho_{i,j}P_{j}(k+1))C_{i}
+\sigma^{2}D_{i}'(\sum_{{j}=1}^L\rho_{i,j}P_{j}(k+1))D_{i}+R_{i}\bigg]u(k),\label{f22}
\end{eqnarray}
i.e.,
\begin{eqnarray}
u(k)=-\Upsilon_{i}(k)^{\dagger}M_{i}(k)x(k).\label{f23}
\end{eqnarray}

From (\ref{x1}),  (\ref{ff6}) and (\ref{f23}),  $\lambda_{k-1}$  can be derived as that
\begin{eqnarray}
\lambda_{k-1}&=&\mbox{E}[(A_{\theta(k)}+B_{\theta(k)}\omega(k))'\lambda_k+Q_{\theta(k)}x(k)|{\cal{G}}_{k-1}]\nonumber\\
&=&\mbox{E}[(A_{\theta(k)}+B_{\theta(k)}\omega(k))'(\sum_{j=1}^L\rho_{i,j}P_{j}(k+1))x(k+1)]+Q_{\theta(k)}x(k)|{\cal{G}}_{k-1}]\nonumber\\
&=&\mbox{E}\bigg[A_{i}'(\sum_{j=1}^L\rho_{i,j}P_{j}(k+1))A_{i}+B_{i}'(\sum_{j=1}^L\rho_{i,j}P_{j}(k+1))B_{i}+Q_{i}\nonumber\\
&&-M_{i}(k)'\Upsilon_{i}(k)^{\dagger}M_{i}(k)|{\cal{G}}_{k-1}\bigg]{x}(k)\nonumber\\
&=&(\sum_{i=1}^L\rho_{s,i}P_{i}(k))x(k),s=1,\cdots,L.\label{f24}
\end{eqnarray}
The  proof  about necessity is end.

(Sufficiency): When the GDRE-MJ (\ref{f7}) has a solution, we will show that Problem$^{\ast}$ is solvable.

Denote $V_{N}(k,x(k))\triangleq\mbox{E}[x(k)'P_{\theta(k)}(k)x(k)]$. From (\ref{x1}) we deduce that
\begin{eqnarray*}
&&V_{N}(k,x(k))-V_{N}(k+1,x(k+1))\\
&=&\mbox{E}[x(k)'P_{\theta(k)}(k)x(k)-x(k+1)'P_{\theta(k+1)}(k+1)x(k+1)]\\
&=&\mbox{E}\bigg\{x(k)'[P_{i}(k)-A_{i}'(\sum_{j=1}^L\rho_{i,j}P_{j}(k+1))A_{i}
-\sigma^{2}B_{i}'(\sum_{{j}=1}^L\rho_{i,j}P_{j}(k+1))B_{i}]x(k)\\
&&-x(k)'[A_{i}'(\sum_{j=1}^L\rho_{i,j}P_{j}(k+1))C_{i}
+\sigma^{2}B_{i}'(\sum_{{j}=1}^L\rho_{i,j}P_{j}(k+1))D_{i}]u(k)\\
&&u(k)'[C_{i}'(\sum_{j=1}^L\rho_{i,j}P_{j}(k+1))A_{i}
+\sigma^{2}D_{i}'(\sum_{{j}=1}^L\rho_{i,j}P_{j}(k+1))B_{i}]x(k)\\
&&-u(k)'[C_{i}'(\sum_{{j}=1}^L\rho_{i,j}P_{j}(k+1))C_{i}
+\sigma^{2}D_{i}'(\sum_{{j}=1}^L\rho_{i,j}P_{j}(k+1))D_{i}]u(k)\bigg\}\\
&=&\mbox{E}\bigg\{x(k)'[Q_{i}-M_{i}(k)'\Upsilon_{i}(k)^{\dagger}M_{i}(k)]x(k)-x(k)'M_{i}'(k)u(k)\\
&&-u(k)'M_{i}(k)x(k)-u(k)'\Upsilon_{i}(k)u(k)+u(k)'R_{i}u(k)\bigg\}\\
&=&\mbox{E}\bigg\{x(k)'Q_{i}x(k)+u(k)'R_{i}u(k)-[u(k)+\Upsilon_{i}(k)^{\dagger}M_{i}(k)x(k)]'
\Upsilon_{i}(k)[u(k)+\Upsilon_{i}(k)^{\dagger}M_{i}(k)x(k)]\bigg\}
\end{eqnarray*}
Adding from $k=0$ to $k=N$ on both sides of the above equation, we have that
\begin{eqnarray}
&&V_{N}(0,x(0))-V_{N}(N+1,x(N+1))\nonumber\\
&=&\mbox{E}\sum_{k=0}^{N}\bigg\{x(k)'Q_{i}x(k)+u(k)'R_{i}u(k)\nonumber\\
&&-[u(k)+\Upsilon_{i}(k)^{\dagger}M_{i}(k)x(k)]'\Upsilon_{i}(k)[u(k)+\Upsilon_{i}(k)^{\dagger}M_{i}(k)x(k)]\bigg\}.\label{f25}
\end{eqnarray}
The above mentioned equation implies that
\begin{eqnarray*}
J_{N}=\mbox{E}[x_{0}'P_{\theta(0)}x_{0}]+\sum_{k=0}^{N}[u(k)+\Upsilon_{i}(k)^{\dagger}M_{i}(k)x(k)]'\Upsilon_{i}(k)[u(k)+\Upsilon_{i}(k)^{\dagger}M_{i}(k)x(k)].
\end{eqnarray*}
Considering $\Upsilon_{i}(k)\geq 0$, we have $J_{N}\geq\mbox{E}[x_{0}'P_{\theta(0)}x_{0}]$. Therefore, the optimal controller can be given by $u(k)=-\Upsilon_{i}(k)^{\dagger}M_{i}(k)x(k)$ and the optimal cost is given by $J_{N}=\mbox{E}[x_{0}'P_{\theta(0)}x_{0}]$.

This completes the proof.

{\bf Remark 2 }\ \ The key technique adopted in this paper is the solving of the FBSDEs-MJ, which is new to our best knowledge. It plays an important role in the design of the optimal controller and stabilization analysis in next section.

\section{Main Result}

The  quadratic optimal control and stabilization problems in infinite horizon will be analyzed in this section. Some necessary definitions will be introduced firstly.

{\bf Definition 3} \ \
The linear system (\ref{x1}) with $u(k)=0$ is asymptotically mean square stable (MSS) if for any initial condition $x_0$ and $\theta(0)$, there holds
\begin{eqnarray}
\lim_{k\rightarrow \infty}\mbox{E}(x(k)'x(k))=0.\nonumber
\end{eqnarray}

{\bf Definition 4} \ \
The system (\ref{x1}) is mean square stabilizable if there is a ${\cal{G}}_k$-measurable controller $u(k)=F_{\theta(k)}x(k)$ satisfying $\lim_{k\rightarrow \infty}\mbox{E}[u(k)'u(k)]=0$, such that system (\ref{x1}) is asymptotically mean square stable.

Denote $A=(A_1,\cdots,A_L), B=(B_1,\cdots,B_L)$. For brevity, we usually say that the pair $(A, B)$ is mean square stabilizable if system (\ref{x1}) is mean square stabilizable.

Now we define the following generalized algebraic Riccati equation with Markov jump as
\begin{eqnarray}
\left\{
\begin{array}{lll}
P_{i}=A_{i}'(\sum_{j=1}^L\rho_{i,j}P_{j})A_{i}
+\sigma^{2}B_{i}'(\sum_{{j}=1}^L\rho_{i,j}P_{j})B_{i}+Q_{i}
-M_{i}'\Upsilon_{i}^{\dagger}M_{i},\\
\Upsilon_{i}\Upsilon_{i}^{\dagger}M_{i}-M_{i}=0,\\
\Upsilon_{i}\geq0,\label{f26}
\end{array}
\right.
\end{eqnarray}
in which
\begin{eqnarray}
\Upsilon_{i}&=&C_{i}'(\sum_{{j}=1}^L\rho_{i,j}P_{j})C_{i}
+\sigma^{2}D_{i}'(\sum_{{j}=1}^L\rho_{i,j}P_{j})D_{i}+R_{i},\label{f27}\\
M_{i}&=&C_{i}'(\sum_{j=1}^L\rho_{i,j}P_{j})A_{i}
+\sigma^{2}D_{i}'(\sum_{{j}=1}^L\rho_{i,j}P_{j})B_{i}.\label{f28}
\end{eqnarray}

For the sake of illustrating the main result, we need to consider the following set which involves linear matric inequality  condition and kernel restriction, whose definition is inspired by \cite{21},
\begin{eqnarray*}
\mathcal{S}\triangleq
\Bigg\{
\begin{smallmatrix}
\tilde{P}=\tilde{P}'\Bigg|
&&\left[
  \begin{smallmatrix}
    A_{i}'(\sum_{j=1}^L\rho_{i,j}\tilde{P}_{j})A_{i}
+\sigma^{2}B_{i}'(\sum_{{j}=1}^L\rho_{i,j}\tilde{P}_{j})B_{i}+Q_{i}-\tilde{P}_{i} & A_{i}'(\sum_{j=1}^L\rho_{i,j}\tilde{P}_{j})C_{i}
+\sigma^{2}B_{i}'(\sum_{{j}=1}^L\rho_{i,j}\tilde{P}_{j})D_{i}\\
    C_{i}'(\sum_{j=1}^L\rho_{i,j}\tilde{P}_{j})A_{i}
+\sigma^{2}D_{i}'(\sum_{{j}=1}^L\rho_{i,j}\tilde{P}_{j})B_{i} & C_{i}'(\sum_{{j}=1}^L\rho_{i,j}\tilde{P}_{j})C_{i}
+\sigma^{2}D_{i}'(\sum_{{j}=1}^L\rho_{i,j}\tilde{P}_{j})D_{i}+R_{i}\\  
  \end{smallmatrix}
\right]\geq 0,\\
&&\mathbf{Ker}(C_{i}'(\sum_{{j}=1}^L\rho_{i,j}\tilde{P}_{j})C_{i}
+\sigma^{2}D_{i}'(\sum_{{j}=1}^L\rho_{i,j}\tilde{P}_{j})D_{i}+R_{i})\subseteq (\mathbf{Ker} C_{i} \cap \mathbf{Ker} D_{i})
\end{smallmatrix}\Bigg\},
\end{eqnarray*}\\
where $\tilde{P}=(\tilde{P}_{1}, \cdots, \tilde{P}_{L})$.

To simplify notation in the sequel, for any $\tilde{P}\in\mathcal{S}$, we denote
\begin{eqnarray}
\left\{
\begin{array}{lll}
\tilde{Q}_{i}=A_{i}'(\sum_{j=1}^L\lambda_{i,j}\tilde{P}_{j})A_{i}
+\sigma^{2}B_{i}'(\sum_{{j}=1}^L\lambda_{i,j}\tilde{P}_{j})B_{i}+Q_{i}-\tilde{P}_{i},\\
\tilde{L}_{i}= A_{i}'(\sum_{j=1}^L\lambda_{i,j}\tilde{P}_{j})C_{i}
+\sigma^{2}B_{i}'(\sum_{{j}=1}^L\lambda_{i,j}\tilde{P}_{j})D_{i},\\
\tilde{R}_{i}= C_{i}'(\sum_{{j}=1}^L\lambda_{i,j}\tilde{P}_{j})C_{i}
+\sigma^{2}D_{i}'(\sum_{{j}=1}^L\lambda_{i,j}\tilde{P}_{j})D_{i}+R_{i}.\label{f29}
\end{array}
\right.
\end{eqnarray}
{\bf Remark 3} \ \  Obviously, we have $\left[
  \begin{array}{cc}
    \tilde{Q}_{i}& \tilde{L}_{i}'\\
    \tilde{L}_{i} & \tilde{R}_{i}\\  
  \end{array}
\right]\geq 0$. In view of Lemma 3, it yields that
\begin{eqnarray} \tilde{R}_{i}\geq 0,\quad  \tilde{Q}_{i}-\tilde{L}_{i}'\tilde{R}_{i}^{\dagger}\tilde{L}_{i}\geq 0, \quad \tilde{L}_{i}'(I-\tilde{R}_{i}\tilde{R}_{i}^{\dagger})=0.\label{f30}
\end{eqnarray}

{\bf Definition 5} \  \
A solution to the GARE-MJ (\ref{f9})-(\ref{f11}) is called a maximal solution, denoted by $P_{max}$, if
\begin{eqnarray}
P_{max}\geq \tilde{P}, \forall \tilde{P}\in {\cal{S}},\label{f31}
\end{eqnarray}
where $P_{max}=(P_{max_1},\cdots,P_{max_L})$.

To make the time horizon $N$ explicit in the finite-horizon LQR problem, we rewrite $\Upsilon_{i}(k), P_i(k)$, and $M_i(k) (i=1,\cdots,L)$ in (\ref{f7})-(\ref{f9}) as $\Upsilon_{i}^N(k), P_i^N(k)$, and $M_i^N(k) (i=1,\cdots,L)$. To facilitate our discussion in the sequels, the terminal weight matrix $P_{j}(N+1)=\tilde{P}_{j}\in \mathcal{S}, j=1,\cdots,L$.

{\bf Definition 6} \ \
Consider the following MJLS with multiplicative noises
\begin{eqnarray}
\left\{
   \begin{aligned}
   x(k+1)&=A_{\theta(k)}x(k)+\omega(k)B_{\theta(k)}x(k),\\
y(k)&=\tilde{Q}_{\theta(k)}^{\frac{1}{2}}x(k),\label{f32}\\
   \end{aligned}
   \right.
\end{eqnarray}
$(A,B,\tilde{Q}^{\frac{1}{2}})$ is said to be exactly observable, if for any $N\geq 0$,
\begin{eqnarray}
y(k)=0, a.s., \forall k\in [0,N]\Rightarrow x_0=0,\nonumber
\end{eqnarray}
where $A=(A_1,\cdots,A_L), B=(B_1,\cdots,B_L), \tilde{Q}^{\frac{1}{2}}=(\tilde{Q}_1^{\frac{1}{2}},\cdots,\tilde{Q}_L^{\frac{1}{2}}).$

{\bf Assumption 1} \ \
$(A,B,\tilde{Q}^{\frac{1}{2}})$ is exactly observable, in which $\tilde{Q}=(\tilde{Q}_{1},\cdots,\tilde{Q}_{L})$ defined as in (\ref{f29}).

{\bf Theorem 2} \ \
Under Assumptions 1 and $\mathcal{S}\neq \emptyset$, if the system (\ref{x1}) is mean square stabilizable, we have the following properties:

For any $k\geq 0$, $P_{i}^{N}(k)$ is convergent when $N\rightarrow\infty$, i.e.,
$\lim\limits_{N\rightarrow\infty}P_{i}^{N}(k)=P_{i}$, in which $P_{i}$ satisfies (\ref{f26})-(\ref{f28}), and $P_{i}$ is the maximal solution to the GARE-MJ.

\emph{Proof}.
Define a new generalized Riccati equation with Markov jump (NGDRE-MJ) with $\theta(k)=i$
\begin{eqnarray}
\left\{
\begin{array}{lll}
X_{i}(k)=A_{i}'(\sum_{j=1}^L\rho_{i,j}X_{j}(k+1))A_{i}
+\sigma^{2}B_{i}'(\sum_{{j}=1}^L\rho_{i,j}X_{j}(k+1))B_{i}+\tilde{Q}_{i}\\
-\tilde{M}_{i}(k)'\tilde{\Upsilon}_{i}(k)^{\dagger}\tilde{M}_{i}(k),\\
\tilde{\Upsilon}_{i}(k)\tilde{\Upsilon}_{i}(k)^{\dagger}\tilde{M}_{i}(k)-\tilde{M}_{i}(k)=0,\\
\tilde{\Upsilon}_{i}(k)\geq0,\label{f33}
\end{array}
\right.
\end{eqnarray}
in which
\begin{eqnarray}
\tilde{\Upsilon}_{i}(k)&=&C_{i}'(\sum_{{j}=1}^L\rho_{i,j}X_{j}(k+1))C_{i}
+\sigma^{2}D_{i}'(\sum_{{j}=1}^L\rho_{i,j}X_{j}(k+1))D_{i}+\tilde{R}_{i},\label{f34}\\
\tilde{M}_{i}(k)&=&C_{i}'(\sum_{j=1}^L\rho_{i,j}X_{j}(k+1))A_{i}
+\sigma^{2}D_{i}'(\sum_{{j}=1}^L\rho_{i,j}X_{j}(k+1))B_{i}+\tilde{L}_i,\label{f35}
\end{eqnarray}
 with its terminal values $X_{i}(N+1)=0$ for $i=1,\cdots,L$ and $\tilde{Q}_{i}, \tilde{L}_{i}, \tilde{R}_{i}$ are denoted as in (\ref{f29}). And the corresponding new cost function can be written as
\begin{eqnarray}
\tilde{J}_{N}=\mbox{E}\sum_{k=0}^{N}\left[
  \begin{array}{cc}
x(k)\\
u(k)
  \end{array}
\right]'\left[
  \begin{array}{cc}
\tilde{Q}_{\theta_{k}} &\tilde{L}_{\theta_{k}}'\\
\tilde{L}_{\theta_{k}} & \tilde{R}_{\theta_{k}}
  \end{array}
\right]\left[
  \begin{array}{cc}
x(k)\\
u(k)
  \end{array}
\right].\label{f36}
\end{eqnarray}
It is clear to know that $\tilde{J}_{N}\geq 0$ from Remark 3.

It is easy to see the difference equation
\begin{eqnarray}
\left\{
\begin{array}{lll}
X_{i}(k)=A_{i}'(\sum_{j=1}^L\rho_{i,j}X_{j}(k+1))A_{i}
+\sigma^{2}B_{i}'(\sum_{{j}=1}^L\rho_{i,j}X_{j}(k+1))B_{i}+\tilde{Q}_{i}\\
-\tilde{M}_{i}(k)'\tilde{\Upsilon}_{i}(k)^{\dagger}\tilde{M}_{i}(k),\\
\tilde{\Upsilon}_{i}(k)\geq0,
\end{array}\label{f37}
\right.
\end{eqnarray}
has a solution, in which $\tilde{\Upsilon}_{i}(k),\tilde{M}_{i}(k)$ are defined as in (\ref{f34}) and (\ref{f35}).

Further we will illustrate that the solution $X_{i}(k)$ of the above equation is positive semi-definite. Considering the following formula
\begin{eqnarray}
\tilde{M}_i'(k)\tilde{\Upsilon}_i^{\dagger}(k)\tilde{M}_i(k)
=-\tilde{M}_i'(k)\tilde{F}_i(k)-\tilde{F}_i'(k)\tilde{M}_i(k)
-\tilde{F}_i'(k)\tilde{\Upsilon}_i(k)\tilde{F}_i(k),\nonumber
\end{eqnarray}
in which $\tilde{F}_i(k)=-\tilde{\Upsilon}_i^{\dagger}(k)\tilde{M}_i(k)$, thus (\ref{f37}) can be rewritten as
\begin{eqnarray*}
X_i^{N}(k)&=&\bar{A}_{i}(k)'(\sum_{j=1}^L\rho_{ij}X_j^{N}(k+1))\bar{A}_{i}(k)
+\bar{B}_{i}(k)'(\sum_{j=1}^L\rho_{ij}X_j^{N}(k+1))\bar{B}_{i}(k)+\bar{Q}_i(k),
\end{eqnarray*}
where
 \begin{eqnarray}
\bar{A}_{i}(k)=A_i+C_i\tilde{F}_i^{N}(k), \bar{B}_{i}(k)=B_i+D_i\tilde{F}_i^{N}(k),\nonumber\\ \bar{Q}_i(k)=\tilde{Q}_i+\tilde{L}_i'\tilde{F}_i(k)
+\tilde{F}_i(k)'\tilde{L}_i+\tilde{F}_i(k)'\tilde{R}_i\tilde{F}_i(k).\label{f38}
\end{eqnarray}

By the Schur complementary, and in view of $\tilde{Q}_{i}\geq0, \tilde{R}_{i}\geq0$,  we have
\begin{eqnarray}
\bar{Q}_i(k)&=&\tilde{Q}_i+\tilde{L}_i'\tilde{F}_i(k)
+\tilde{F}_i(k)'\tilde{L}_i+\tilde{F}_i(k)'\tilde{R}_i\tilde{F}_i(k)\nonumber \\&\geq&\tilde{L}_i'\tilde{R}_i^{\dagger}\tilde{L}_i+\tilde{L}_i'\tilde{R}_i^{\dagger}\tilde{R}_i\tilde{F}_i(k)
+\tilde{F}_i(k)'\tilde{R}_i\tilde{R}_i^{\dagger}\tilde{L}_i+\tilde{F}_i(k)'\tilde{R}_i\tilde{R}_i^{\dagger}\tilde{R}_i\tilde{F}_i(k)\nonumber\\
&=&(\tilde{L}_i+\tilde{R}_i\tilde{F}_i(k))'\tilde{R}_i^{\dagger}(\tilde{L}_i+\tilde{R}_i\tilde{F}_i(k))\geq 0,\label{f39}
\end{eqnarray}
and on the ground of $X_{i}^{N}(N+1)=0,i=1,\cdots,L$, it yields that $X_{i}^{N}(N)\geq0$ and by induction, it is not hard to verify that $X_{i}^{N}(k)\geq0$, for $0\leq k \leq N$.

 Next we will investigate $\tilde{\Upsilon}_{i}(k)\tilde{\Upsilon}_{i}(k)^{\dagger}\tilde{M}_{i}(k)-\tilde{M}_{i}(k)=0$.

Considering $\tilde{\Upsilon}_{i}\geq0$, it yields that
\begin{eqnarray}
\tilde{\Upsilon}_{i}^{\dagger}=V_{i}\left[
  \begin{array}{cc}
    T_{i}^{-1} & 0 \\
    0& 0
  \end{array}
\right]V_{i}',\label{f40}
\end{eqnarray}
where $T_{i}>0$  has same dimension with the rank of $\tilde{\Upsilon}_{i}$ and $V_{i}$ is an orthogonal matrix.
Now, let $V_{i}$ decompose as
$\left[
  \begin{array}{cc}
V_{i}^{1} & V_{i}^{2}
  \end{array}
\right]$
where the columns of the matrix $V_{i}^{2}$ form a basis of $\mathbf{Ker}(\tilde{\Upsilon}_{i})$. The positive semi-definite of matrices $\tilde{R}_{i}, \sum_{{j}=1}^L\rho_{i,j}X_{j}(k+1)$ yields that $\mathbf{Ker}(\tilde{\Upsilon}_{i})\subseteq \mathbf{Ker}(\tilde{R}_{i})$. A simple calculation yields that
\begin{eqnarray*}
&&[A_{i}'(\sum_{j=1}^L\rho_{i,j}X_{j}(k+1))C_{i}
+\sigma^{2}B_{i}'(\sum_{{j}=1}^L\rho_{i,j}X_{j}(k+1))D_{i}+\tilde{L}_{i}'][I-\tilde{\Upsilon}_{i}(k)\tilde{\Upsilon}_{i}(k)^{\dagger}]\\
&=&[A_{i}'(\sum_{j=1}^L\rho_{i,j}(X_{j}(k+1)+\tilde{P}_{j}))C_{i}
+\sigma^{2}B_{i}'(\sum_{{j}=1}^L\rho_{i,j}(X_{j}(k+1)+\tilde{P}_{j}))D_{i}]V_{i}^{2}(V_{i}^{2})'.
\end{eqnarray*}
On the ground of $\tilde{\Upsilon}_{i}(k)V_{i}^{2}(V_{i}^{2})'=0$, it is easy to verify $\tilde{R}_{i}V_{i}^{2}(V_{i}^{2})'=0$. And further considering the condition of $\mathbf{Ker}(C_{i}'(\sum_{{j}=1}^L\rho_{i,j}\tilde{P}_{j})C_{i}
+\sigma^{2}D_{i}'(\sum_{{j}=1}^L\rho_{i,j}\tilde{P}_{j})D_{i}+R_{i})\subseteq (\mathbf{Ker} C_{i} \cap \mathbf{Ker} D_{i})$, we have $\tilde{\Upsilon}_{i}(k)\tilde{\Upsilon}_{i}(k)^{\dagger}\tilde{M}_{i}(k)-\tilde{M}_{i}(k)=0$. Up to now, we know that NGDRE-MJ (\ref{f33}) exists a positive semi-definite solution. Considering the result of Theorem 1, the optimal controller and cost value of the new cost function subject to (\ref{x1}) are $u^{\ast}(k)=-\tilde{\Upsilon}_{i}(k)^{\dagger}\tilde{M}_{i}(k)x(k)$ and
$\tilde{J}_{N}^{\ast}=\mbox{E}[x_{0}'X_{\theta(0)}^{N}(0)x_{0}]$, respectively.
 Thus, for any $N$, we have
 \begin{eqnarray*}
\tilde{J}_{N}^{\ast}=\mbox{E}[x_{0}'X_{\theta(0)}^{N}(0)x_{0}]\leq\mbox{E}[x_{0}'X_{\theta(0)}^{N+1}(0)x_{0}]=\tilde{J}_{N+1}^{\ast}.
\end{eqnarray*}
The arbitrariness of $x_{0}$ implies that $X_{\theta(0)}^{N}(0)$ increases with respect to $N$. Next, we will show the boundedness of $X_{\theta(0)}^{N}(0)$. Since system (\ref{x1}) is stabilizable in the mean square sense, there exists $u(k)=F_{\theta(k)}x(k)$ satisfying
\begin{eqnarray*}
\lim\limits_{k\rightarrow\infty}\mbox{E}(x(k)'x(k))=0.
\end{eqnarray*}
Hence, we have that
\begin{eqnarray*}
\tilde{J}_{N}^{\ast}\leq\tilde{J}&=&\mbox{E}\sum_{k=0}^{\infty}\left[
  \begin{array}{cc}
x(k)\\
u(k)
  \end{array}
\right]'\left[
  \begin{array}{cc}
\tilde{Q}_{\theta_{k}} &\tilde{L}_{\theta_{k}}'\\
\tilde{L}_{\theta_{k}} & \tilde{R}_{\theta_{k}}
  \end{array}
\right]\left[
  \begin{array}{cc}
x(k)\\
u(k)
  \end{array}
\right]\\
&=&\mbox{E}\sum_{k=0}^{\infty}\Bigg\{x(k)'\left[
  \begin{array}{cc}
I\\
F_{\theta(k)}
  \end{array}
\right]'\left[
  \begin{array}{cc}
\tilde{Q}_{\theta_{k}} &\tilde{L}_{\theta_{k}}'\\
\tilde{L}_{\theta_{k}} & \tilde{R}_{\theta_{k}}
  \end{array}
\right]\left[
  \begin{array}{cc}
I\\
F_{\theta(k)}
  \end{array}
\right]x(k)\Bigg\}\\
&&\leq\lambda_{max}\mbox{E}(x(k)'x(k))\\
&&\leq\lambda_{max}\cdot c\cdot\mbox{E}(x_{0}'x_{0}),
\end{eqnarray*}
where $\lambda_{max}$ denotes the maximum eigenvalue of
$\left[
  \begin{array}{cc}
I\\
F_{\theta(k)}
  \end{array}
\right]'\left[
  \begin{array}{cc}
\tilde{Q}_{\theta_{k}} &\tilde{L}_{\theta_{k}}'\\
\tilde{L}_{\theta_{k}} & \tilde{R}_{\theta_{k}}
  \end{array}
\right]\left[
  \begin{array}{cc}
I\\
F_{\theta(k)}
  \end{array}
\right]$
and $c$ is a positive constant. The above formula implies that
\begin{eqnarray*}
\mbox{E}[x_{0}'X_{\theta(0)}^{N}(0)x_{0}]\leq\lambda_{max}\cdot c\cdot\mbox{E}(x_{0}'x_{0}),
\end{eqnarray*}
i.e., $X_{\theta(0)}^{N}(0)\leq \lambda_{max}\cdot c$ in view of the arbitrariness of $x_{0}$.

Up to now, we can say that $X_{\theta(0)}^{N}(0)$ is bounded. In considering of the monotonicity of $X_{\theta(0)}^{N}(0)$, we deduce that $X_{\theta(0)}^{N}(0)$ is convergent. Note that the variables given in NGDRE-MJ are time invariant for $N$ due to the choice of $X_{j}^{N+1}=0$, so we have
\begin{eqnarray*}
\lim\limits_{N\rightarrow\infty}X_{i}^{N}(k)=\lim\limits_{N\rightarrow\infty}X_{i}^{N-k}(0)=X_{i}, i=1,\cdots , L,
\end{eqnarray*}
at the same time, we have that
\begin{eqnarray*}
\lim\limits_{N\rightarrow\infty}\tilde{\Upsilon}_{i}^{N}(k)&=&\tilde{\Upsilon}_{i},\\
 \lim\limits_{N\rightarrow\infty}\tilde{M}_{i}^{N}(k)&=&\tilde{M}_{i}, \ \ i=1,\cdots , L.
\end{eqnarray*}
Therefore, we can say that $X_{i}$ is a solution of the following NGARE-MJ
\begin{eqnarray}
\left\{
\begin{array}{lll}
X_{i}=A_{i}'(\sum_{j=1}^L\rho_{i,j}X_{j})A_{i}
+\sigma^{2}B_{i}'(\sum_{{j}=1}^L\rho_{i,j}X_{j})B_{i}+\tilde{Q}_{i}
-\tilde{M}_{i}'\tilde{\Upsilon}_{i}^{\dagger}\tilde{M}_{i},\\
\tilde{\Upsilon}_{i}\tilde{\Upsilon}_{i}^{\dagger}\tilde{M}_{i}-\tilde{M}_{i}=0,\\
\tilde{\Upsilon}_{i}\geq0,
\end{array}\label{f41}
\right.
\end{eqnarray}
in which
\begin{eqnarray}
\tilde{\Upsilon}_{i}&=&C_{i}'(\sum_{{j}=1}^L\rho_{i,j}X_{j})C_{i}
+\sigma^{2}D_{i}'(\sum_{{j}=1}^L\rho_{i,j}X_{j})D_{i}+\tilde{R}_{i},\label{f42}\\
\tilde{M}_{i}&=&C_{i}'(\sum_{j=1}^L\rho_{i,j}X_{j})A_{i}
+\sigma^{2}D_{i}'(\sum_{{j}=1}^L\rho_{i,j}X_{j})B_{i}+\tilde{L}_i.\label{f43}
\end{eqnarray}

Next we will mainly illustrate that $X_{i}>0$.  Owing to the positive semi-definiteness of $X_{i}^{N}(k)$,  its limit $X_{i}$ is also positive semi-definite, i.e., $X_{i}\geq0$.  Now we verify $X_{i}>0$. If not, there must exist nonzero vector $x_{0}$ such that $E[x_{0}'X_{i}x_{0}]=0$.

Define  the Lyapunov function as
\begin{eqnarray}
V_X(k,x(k))=\mbox{E}[x(k)'X_{\theta(k)}x(k)],k\geq 0.\label{f44}
\end{eqnarray}
Since $X_{\theta(k)}=X_i$ for $\theta(k)=i,i=1,\cdots,L$, then $V_X(k,x(k))\geq 0$. So we have
\begin{eqnarray}
&&\sum_{k=0}^{N}[V_X(k+1,x(k+1))-V_X(k,x(k))]\nonumber\\
&=&V_X(N+1,x(N+1))-V_X(0,x_{0})\nonumber\\
&=&\mbox{E}\{x(N+1)'X_{i}x(N+1)-x_{0}'X_{i}x_{0}\}\nonumber\\
&=&-\sum_{k=0}^{N}\mbox{E}\{x(k)'\tilde{Q}_ix(k)+x(k)'\tilde{L}_i'u(k)+u(k)'\tilde{L}_ix(k)+u(k)'\tilde{R}_iu(k)\}\nonumber\\
&=&-\sum_{k=0}^{N}\mbox{E}\{x(k)'\tilde{Q}_ix(k)+x(k)'\tilde{L}_i'\tilde{F}_{i}x(k)
+x(k)'\tilde{F}_{i}'\tilde{L}_ix(k)+x(k)'\tilde{F}_{i}'\tilde{R}_i\tilde{F}_{i}x(k)\}\nonumber\\
&=&-\sum_{k=0}^{N}\mbox{E}\{x(k)'[\tilde{Q}_i+\tilde{L}_i'\tilde{F}_{i}
+\tilde{F}_{i}'\tilde{L}_i+\tilde{F}_{i}'\tilde{R}_i\tilde{F}_{i}]x(k)\}\nonumber\\
&=&-\sum_{k=0}^{N}\mbox{E}\{x(k)'\bar{Q}_{i}x(k)\}\leq0,\label{f45}
\end{eqnarray}
where $u(k)=\tilde{F}_{i}x(k)=-\tilde{\Upsilon}_i^{\dagger}\tilde{M}_ix(k)$ is used in the fourth equation.
Obviously,
\begin{eqnarray*}
0\leq\sum_{k=0}^{N}\mbox{E}\{x(k)'\bar{Q}_{i}x(k)\}=-\mbox{E}[x(N+1)'X_{i}x(N+1)]\leq0,
\end{eqnarray*}
it implies that
\begin{eqnarray*}
\sum_{k=0}^{N}\mbox{E}\{x(k)'\bar{Q}_{i}x(k)\}=0,
\end{eqnarray*}
i.e.,
\begin{eqnarray}
\bar{Q}_{i}^{\frac{1}{2}}x(k)=0.\label{f46}
\end{eqnarray}
 From Theorem 4 and Proposition 1 in \cite{23}, we know that  the  exact observable of $(\bar{A},\bar{B},\bar{Q}^{\frac{1}{2}})$ can be deduced by the exact observable of $(A, B, \tilde{Q}^{\frac{1}{2}})$, in which $\bar{A}=A_i+C_i\tilde{F}_i$, $\bar{B}=B_i+D_i\tilde{F}_i$. Therefore, from (\ref{f46}), it yields that $x_{0}=0$ which is contrary with $x_{0}\neq0$. Hence, we have $X_{i}>0$.

Define $P_{i}^{N}(k)=X_{i}^{N}(k)+\tilde{P}_{i}$. It is easy to verify that $P_{i}^{N}(k)$ satisfies the GDRE-MJ (\ref{f7}) and monotonically increasing with respect to $N$ and bounded. Therefore, there exists a constant $P_{i}$ satisfying
\begin{eqnarray*}
P_{i}=\lim\limits_{N\rightarrow\infty}X_{i}^{N}(k)+\tilde{P}_{i}=X_{i}+\tilde{P}_{i}.
\end{eqnarray*}
Obviously, $P_{i}$ satisfies GARE-MJ (\ref{f26}). Moreover, for the arbitrariness of $\tilde{P}_{i}$ and $X_{i}>0$, we can obtain that $P_{i}\geq \tilde{P}_{i}$, i.e., $P_{i}$ is the maximal solution to the GARE-MJ (\ref{f26}). The proof is complete.

{\bf Remark 4} \ \
The above proof implies that the solvability of the GARE-MJ (\ref{f26}) is equivalent to the solvability of the NGARE-MJ (\ref{f41}).


{\bf Theorem 3} \ \
If ${\cal{S}}\neq \emptyset$ and Assumption 1 are satisfied, then the closed-loop system (\ref{x1}) is mean-square stabilizable if and only if the GARE-MJ (\ref{f26}) has a solution $P=(P_1,\cdots,P_L)$, which is also the maximal solution to the GARE-MJ (\ref{f26}).

In this case, the optimal stabilizing solution is given by
\begin{eqnarray}
u^{\ast}(k)=F_{\theta(k)}x(k),\label{f47}
\end{eqnarray}
where $F_{\theta(k)}=F_i$ for $\theta(k)=i, i=1,\cdots,L$, and
\begin{eqnarray}
F_i&=&-[C_{i}'(\sum_{{j}=1}^L\rho_{i,j}P_{j})C_{i}
+\sigma^{2}D_{i}'(\sum_{{j}=1}^L\rho_{i,j}P_{j})D_{i}+R_{i}]^{\dagger}\nonumber\\
&&\times [C_{i}'(\sum_{j=1}^L\rho_{i,j}P_{j})A_{i}
+\sigma^{2}D_{i}'(\sum_{{j}=1}^L\rho_{i,j}P_{j})B_{i}],\label{f48}
\end{eqnarray}
and the optimal cost functional is
\begin{eqnarray}
J^{*}=\mbox{E}[x_0'P_{\theta(0)}x_0].\label{f49}
\end{eqnarray}

\emph{Proof}. "Sufficiency": We will show that under the conditions of $\cal{S} \neq \emptyset$ and Assumption 1, when the GARE-MJ (\ref{f26}) has a solution $P$, the closed-loop system (\ref{x1}) is mean-square stabilizable.

Let $\tilde{P}\in {\cal{S}}$. From Remark 4, when the GARE-MJ (\ref{f26}) has a solution $P$, the NGARE-MJ (\ref{f41}) has a positive definite solution $X=(X_1,\cdots, X_L)$, i.e., $X_i>0, i=1,\cdots,L$. Moreover, $P=X+\tilde{P}$. Next, we will show that system (\ref{x1}) with $u(k)=F_{\theta(k)}x(k)$, where $F_{\theta(k)}=F_i,i=1,\cdots,L$, is denoted by (\ref{f48}), i.e.,
\begin{eqnarray}
x(k+1)=(A_{\theta(k)}+C_{\theta(k)}F_{\theta(k)})x(k)+\omega(k)(B_{\theta(k)}+D_{\theta(k)}F_{\theta(k)})x(k)\label{f50}
\end{eqnarray}
is mean square stabilizable. In view of the relationship $F_i=\tilde{F}_i(i=1,\cdots,L)$, we can see that the stabilization for the system (\ref{x1}) with $u(k)=F_{\theta(k)}x(k)$ is equivalent to the stabilization for the system (\ref{x1}) with $u(k)=\tilde{F}_{\theta(k)}x(k)$. We define the Lyapunov function as
\begin{eqnarray}
V_X(k,x(k))=\mbox{E}[x(k)'X_{\theta(k)}x(k)],k\geq 0.\label{f51}
\end{eqnarray}
Since $X_{\theta(k)}=X_i$ for $\theta(k)=i,i=1,\cdots,L$, then $V_X(k,x(k))\geq 0$. So we have
\begin{eqnarray}
&&V_X(k+1,x(k+1))-V_X(k,x(k))\nonumber\\
&=&\mbox{E}\{x(k+1)'X_{\theta(k+1)}x(k+1)-x(k)'X_{\theta(k)}x(k)\}\nonumber\\
&=&-\mbox{E}\{x(k)'\tilde{Q}_ix(k)+x(k)'\tilde{L}_i'u(k)+u(k)'\tilde{L}_ix(k)+u(k)'\tilde{R}_iu(k)\}\nonumber\\
&=&-\mbox{E}\Bigg\{\left[
                    \begin{array}{c}
                      x(k) \\
                      u(k) \\
                    \end{array}
                  \right]'
\left[
                            \begin{array}{cc}
                              \tilde{Q}_i & \tilde{L}_i' \\
                              \tilde{L}_i & \tilde{R}_i \\
                            \end{array}
                          \right]\left[
                    \begin{array}{c}
                      x(k) \\
                      u(k) \\
                    \end{array}
                  \right]
\Bigg\}\label{f52}\\
&\leq & 0,\nonumber
\end{eqnarray}
where $\tilde{Q}_i,\tilde{L}_i,\tilde{R}_i$ defined in (\ref{f29}).

From (\ref{f52}), we can see that $V_{X}(k,x(k))$ is non-increasing with respect to $k$. That implies $V_X(k,x(k))\leq V_X(0,x_0)$, i.e., $V_X(k,x(k))$ is bounded. Therefore, $\lim_{k\rightarrow +\infty}V_X(k,x(k))$ exists.

Now for any integer $m\geq0$, taking  summation from $k=m$ to $k=m+N$ on both sides of the upper formulation, we can obtain that
\begin{eqnarray}
&&V_X(m+N+1,x(m+N+1))-V_X(m,x(m))\nonumber\\
&=&-\sum_{k=m}^{m+N}\mbox{E}\Bigg\{\left[
                    \begin{array}{c}
                      x(k) \\
                      u(k) \\
                    \end{array}
                  \right]'
\left[
                            \begin{array}{cc}
                              \tilde{Q}_i & \tilde{L}_i' \\
                              \tilde{L}_i & \tilde{R}_i \\
                            \end{array}
                          \right]\left[
                    \begin{array}{c}
                      x(k) \\
                      u(k) \\
                    \end{array}
                  \right]
\Bigg\}.\label{f53}
\end{eqnarray}
In view of the convergence of $V_X(k,x(k))$, when we take limitation of $m$ on both sides of the aforementioned equation, the following result can be derived as
\begin{eqnarray}
&&-\lim\limits_{m\rightarrow\infty}\sum_{k=m}^{m+N}\mbox{E}\Bigg\{\left[
                    \begin{array}{c}
                      x(k) \\
                      u(k) \\
                    \end{array}
                  \right]'
\left[
                            \begin{array}{cc}
                              \tilde{Q}_i & \tilde{L}_i' \\
                              \tilde{L}_i & \tilde{R}_i \\
                            \end{array}
                          \right]\left[
                    \begin{array}{c}
                      x(k) \\
                      u(k) \\
                    \end{array}
                  \right]
\Bigg\}\nonumber\\
&=&\lim\limits_{m\rightarrow\infty}\left[V_1(m+N+1,x(m+N+1))-V_1(m,x(m))\right]\nonumber\\
&=&0.\label{f54}
\end{eqnarray}

Further considering that the optimal cost function of $\tilde{J}_{N}$ is $\mbox{E}[x_{0}'X_{\theta(0)}^{N}x_{0}]\geq 0$, via a time-shift of length of $m$, it yields that
\begin{eqnarray}
0&=&\lim\limits_{m\rightarrow\infty}\sum_{k=m}^{m+N}\mbox{E}\Bigg\{\left[
                    \begin{array}{c}
                      x(k) \\
                      u(k) \\
                    \end{array}
                  \right]'
\left[
                            \begin{array}{cc}
                              \tilde{Q}_i & \tilde{L}_i' \\
                              \tilde{L}_i & \tilde{R}_i \\
                            \end{array}
                          \right]\left[
                    \begin{array}{c}
                      x(k) \\
                      u(k) \\
                    \end{array}
                  \right]
\Bigg\}\nonumber\\
&&\geq \lim\limits_{m\rightarrow\infty}\mbox{E}[x_{m}'X_{\theta(m)}^{m+N}x_{m}]\nonumber\\
&&=\lim\limits_{m\rightarrow\infty}\mbox{E}[x_{m}'X_{\theta(0)}^{N}x_{m}]\nonumber\\
&&\geq 0.\label{f55}
\end{eqnarray}

Obviously, it implies that $\lim\limits_{m\rightarrow\infty}\mbox{E}[x_{m}'X_{\theta(0)}^{N}x_{m}]=0$. On the ground of the positive definiteness of $X_{\theta(0)}^{N}$, it is easy to verify that $\lim\limits_{m\rightarrow\infty}\mbox{E}[x_{m}'x_{m}]=0$. That is to say that the controller $u(k)=\tilde{F}_{\theta(k)}x(k)=F_{\theta(k)}x(k)$ stabilizes system (\ref{x1}) in the mean square sense.

Lastly, we show the optimal controller and optimal cost. Define
\begin{eqnarray}
V_P(k,x(k))=\mbox{E}[x(k)'P_{\theta(k)}x(k)],k\geq 0.\label{f56}
\end{eqnarray}
 Therefore,
\begin{eqnarray}
&&\sum_{k=0}^{N}[V_P(k+1,x(k+1))-V_P(k,x(k))]\nonumber\\
&=&\sum_{k=0}^{N}\mbox{E}[(u(k)+\Upsilon_{\theta(k)}(k)^{\dagger}M_{\theta(k)}(k)x(k))'\Upsilon_{\theta(k)}(k)
(u(k)+\Upsilon_{\theta(k)}(k)^{\dagger}M_{\theta(k)}(k)x(k))]\nonumber\\
&&-\sum_{k=0}^{N}\mbox{E}[x(k)'Q_{\theta(k)}x(k)+u(k)'R_{\theta(k)}u(k)],\label{f57}
\end{eqnarray}
that is,
\begin{eqnarray}
J_{N}&=&\mbox{E}[x_{0}'P_{\theta(0)}x_{0}]+\sum_{k=0}^{N}\mbox{E}[(u(k)+\Upsilon_{\theta(k)}(k)^{\dagger}M_{\theta(k)}(k)x(k))'\Upsilon_{\theta(k)}(k)\nonumber\\
&&\cdot(u(k)+\Upsilon_{\theta(k)}(k)^{\dagger}M_{\theta(k)}(k)x(k))],\label{f58}
\end{eqnarray}
where $\Upsilon_{\theta(k)}(k)$ and $M_{\theta(k)}(k)$ are as in (\ref{f8}) and (\ref{f9}).

Thus the infinite cost function can be denoted as the following form on account of the mean square stabilizable of system (\ref{x1}), i.e.,
\begin{eqnarray}
J=\mbox{E}[x_{0}'P_{\theta(0)}x_{0}]+\sum_{k=0}^{\infty}\mbox{E}[(u(k)+\Upsilon_{\theta(k)}^{\dagger}M_{\theta(k)}x(k))'\Upsilon_{\theta(k)}
(u(k)+\Upsilon_{\theta(k)}^{\dagger}M_{\theta(k)}x(k))].\label{f59}
\end{eqnarray}
Hence, it's tempting to conclude that the optimal controller is $u^{\ast}(k)=-\Upsilon_{\theta(k)}^{\dagger}M_{\theta(k)}x(k)$ and furthermore the corresponding optimal cost is $J^{\ast}=\mbox{E}[x_{0}'P_{\theta(0)}x_{0}]$.

``Necessity": On the other hand, we should illustrate that if ${\cal{S}}\neq \emptyset$, when the closed-loop system (\ref{x1}) is mean square stabilizable, then the GARE-MJ (\ref{f26}) has a solution $P=(P_1,\cdots,P_L)$, which is also the maximal solution to the GARE-MJ (\ref{f26}). In fact, the existence of solutions to the GARE-MJ (\ref{f26}) have been obtained in Theorem 2. The proof is complete.

{\bf Remark 5} \ \  In \cite{14}, their conclusions can be summarized that under the precondition that the system is stabilizable, based on some positive semi-definite and kernel restrictions on some matrices, necessary and sufficient conditions about the existence of the mean
square stabilizing solution  for a set of
generalized coupled algebraic Riccati equations (GCARE) is derived, i.e., this conclusion was only to study the existence of the stabilizing solution of the GCARE on the basis of the stabilization,  moreover, the given conditions are all operator type which is not easy to be tested.  However, compared with it,  the result expressed by Theorem 3 is clearly illustrated the  stabilization problem with indefinite weighting matrices by the method of transformation that  the stabilization problem of indefinite case is reduced to a definite one whose stabilization condition is expressed by defining Lyapunov function via the optimal cost subject to a new algebraic Riccati equation involving Markov jump.

\section{Numerical example}

A numerical example will be given in this section to further illustrate our result. Now considering system (\ref{x1}) with the following coefficients as
\begin{eqnarray*}
&&A_{1}=\frac{1}{2}, B_{1}=-\frac{1}{2}, C_{1}=\frac{1}{2}, D_{1}=-\frac{1}{2}, Q_{1}=-1, R_{1}=-3;\\
&&A_{2}=\frac{1}{4}, B_{2}=-\frac{1}{4}, C_{2}=\frac{1}{4}, D_{2}=-\frac{1}{4}, Q_{2}=20, R_{2}=0.
\end{eqnarray*}
The transition probabilities of the Markov chain $\{\theta(k); k=1,2,\cdots\}$ taking value in $\{1,2\}$ are $\lambda_{11}=0.2$ and $\lambda_{22}=0.6$. The variance of system noise is 1.

By simply computing we know that
\begin{eqnarray*}
\tilde{Q}_{1}&=&-0.9\tilde{P}_{1}+0.4\tilde{P}_{2}-1,\ \ \ \tilde{L}_{1}=0.1\tilde{P}_{1}+0.4\tilde{P}_{2},\ \ \ \tilde{R}_{1}=0.1\tilde{P}_{1}+0.4\tilde{P}_{2}-3,\\
\tilde{Q}_{2}&=&0.05\tilde{P}_{1}-0.925\tilde{P}_{2}+20,\ \ \ \tilde{L}_{2}=0.05\tilde{P}_{1}+0.075\tilde{P}_{2},\ \ \ \tilde{R}_{2}=0.05\tilde{P}_{1}+0.075\tilde{P}_{2}.
\end{eqnarray*}
Firstly, we calculate the solution of set $\mathcal{S}$. In the case of $\tilde{R}_{2}=0$, i.e., $\tilde{P}_{1}=-\frac{3}{2}\tilde{P}_{2}$, it is easy to verify that the condition of
$\mathbf{Ker}(\tilde{R}_{2})\subseteq (\mathbf{Ker} C_{2} \cap \mathbf{Ker} D_{2})$ is not satisfied. Hence, $\tilde{R}_{2}\neq0$. In this situation, by Schurs Lemma, the matric inequality in set $\mathcal{S}$ can be equivalent to write as
\begin{eqnarray}
\left\{
\begin{array}{lll}
-0.1\tilde{P}_{1}^{2}+2.6\tilde{P}_{1}-1.6\tilde{P}_{2}-0.4P_{1}\tilde{P}_{2}+3\geq 0\\
-0.075\tilde{P}_{2}^{2}+\tilde{P}_{1}+1.5\tilde{P}_{2}-0.05\tilde{P}_{1}\tilde{P}_{2}\geq 0\\
0.1\tilde{P}_{1}+0.4\tilde{P}_{2}-3>0\\
0.05\tilde{P}_{1}+0.075\tilde{P}_{2}>0\\
-0.9\tilde{P}_{1}+0.4\tilde{P}_{2}-1\geq0\\
0.05\tilde{P}_{1}-0.925\tilde{P}_{2}+20\geq0,
\end{array}
\right.
\end{eqnarray}
\begin{figure}[htbp]
  \begin{center}
  \includegraphics[width=1\textwidth]{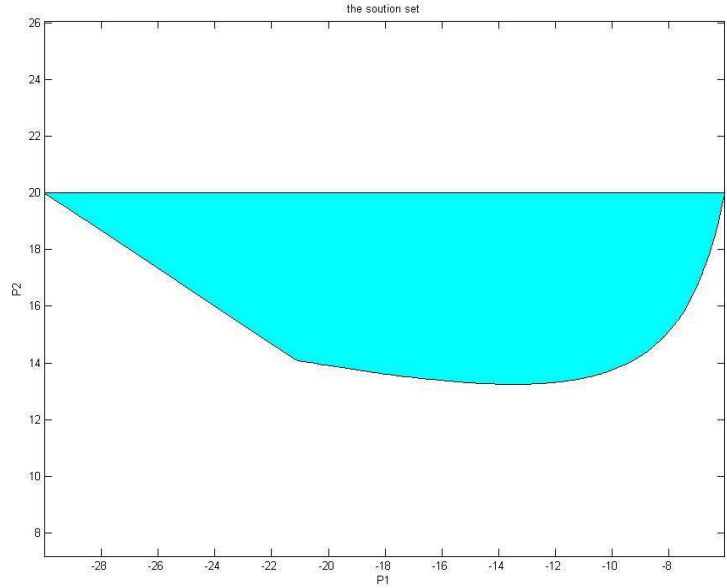}
  \caption{the solution set.} \label{fig:digit}
  \end{center}
\end{figure}
 we can obtain the solution set of the above inequalities be shown as in Fig. 1 by using MATLAB tool. Further it is obviously to see that for any $\tilde{P}=(\tilde{P}_{1},\tilde{P}_{2})$ in above solution set, the condition of $\mathbf{Ker}(\tilde{R}_{i})\subseteq (\mathbf{Ker} C_{i} \cap \mathbf{Ker} D_{i})$, $i=1,2$ is satisfied. Therefore, $\mathcal{S}\neq\emptyset$ and the solution set is expressed as in Fig. 1. And the condition of Assumption 1 can be easy tested. Further the GARE-MJ (\ref{f26}) can be solved as $P_{max}=(\sqrt{103}-11, 20)$. Therefore, $F_{1}=-1.61$ and $F_{2}=-1$. When $\theta(k)=2$, that is, the optimal controller is $u(k)=-x(k)$, in this case, $x(k)=0$, obviously the system is stabilized in the mean square sense. On condition of $\theta(k)=1$, the optimal controller can be given as $u(k)=-1.61x(k)$ and the simulation result is shown in Fig. 2. It can be seen
that the state $x(k)$ is stabilized with the optimal controller, as expected.

 \begin{figure}[htbp]
  \begin{center}
  \includegraphics[width=1\textwidth]{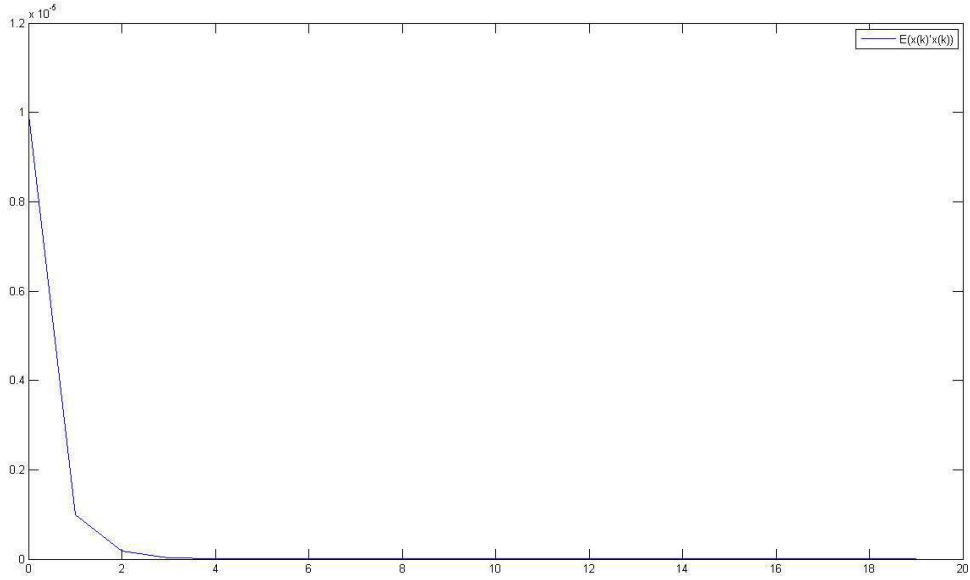}
  \caption{Simulations for the state trajectory $E[x'(k)x(k)]$.} \label{fig:digit}
  \end{center}
\end{figure}

\section{Conclusions}         
This article mainly study the linear quadratic optimal control and
 stabilization problem for discrete-time systems involving Markov jump and multiplicative noise. The state and control weighting matrices in the cost function are allowed to be indefinite.
By solving the FBSDEs-MJ derived from the extended maximum principle, we conclude that the indefinite optimal control problem in finite-horizon is solvable if and only if the corresponding GDRE-MJ has a solution, which is an easy verifiable conclusion compared with operator type results. What's more, in this article we first develop the necessary and sufficient conditions that stabilize the Markov jump discrete-time systems in the mean square sense with indefinite weighting matrices in the cost. More concretely, based only on linear matrix inequality and kernel restrictions, under the basic assumption that the system is exactly observable, the stabilization of Markov jump system can be equivalent to the existence of the maximum solution for the GARE-MJ.
Finally, we give a numerical example to illustrate our correctness of the main result.

\bibliographystyle{plain}        

\begin{thebibliography}{99}     


\bibitem{1}  O. L. V. Costa, M. D. Fragoso and R. P. Marques, ``Discrete Time Markov Jump Linear
Systems," New York: Springer-Verlag, 2005.

\bibitem{2}  J. K. Tugnait, ``Adaptive estimation and identification for discrete
systems with Markov jump parameters," IEEE Trans. Autom. Control,
vol. 27. pp. 1054-1064, 1982.

\bibitem{3}   P. Shi, E. K. Boukas, and R. K. Agarwal, ``Control of Markovian jump
discrete-time systems with norm bounded uncertainty and unknown
delay," IEEE Trans. Autom. Control, vol. 44, no. 11, pp. 2139-2144,
 1999.

\bibitem{4} L. Wu, P. Shi, and H. Gao, ``State estimation and sliding mode control
of Markovian jump singular systems," IEEE Trans. Autom. Control,
vol. 55, no. 5, pp. 1213-1219,  2010.

\bibitem{5}  O. L. V. Costa, E. O. Assumpia Filho, E. K. Boukas, and R. P. Marques, ``Constrained
quadratic state feedback control of discrete-time Markovian jump linear systems," Automatica, vol. 35, pp. 617-626, 1999.

\bibitem{6}   O. L. V. Costa, ``Linear minimum mean square error estimation for discrete-time
Markovian jump linear systems," IEEE Trans. Autom. Control,
vol. 39, no. 8, pp. 1685-1689, 1994.

\bibitem{7} O. L. V. Costa and J. B. R. do Val, ``Full information H control for
discrete-time infinite Markovjump parameter systems," J. Math. Anal.
Appl., vol. 202, pp. 578-603, 1996.

\bibitem{8}  O. L. V. Costa,  M. D. Fragoso, ``Discrete-time LQ-optimal control problems for infnite
Markov jump parameter systems,"  IEEE  Trans. Autom. Control, vol. 40, no. 12, pp. 2076-2088, 1995.

\bibitem{9} J. B. R. do Val and E. F. Costa. ``Stabilizability and positiveness
of solutions of the jump linear quadratic problem and the coupled
algebraic Riccati equation," Trans. Autom. Control, vol. 50, no. 5, pp. 691-695, 2005.

\bibitem{10}  S. Chen, X. Li, and X. Zhou, ``Stochastic linear quadratic
regulators with indefinite control weight costs," SIAM J.
Control Optim., vol.36, no.5, pp. 1685-1702, 1998.

\bibitem{11}  X. Li, and X. Zhou, ``Indefinite stochastic LQ controls with Markovian jumps in a finite time horizon,"
Commun. Inf. Syst. vol.2, no.3, pp., 265-282,  2002.

\bibitem{12}  X. Li, X. Zhou, and M. A. Rami, ``Indefinite stochastic
linear quadratic control with Markovian jumps in infinite
time horizon," Journal of Global Optimization., vol.27, pp.
149-175, 2003.

\bibitem{13} O. L. V. Costa and L. P. Wanderlei, ``Indefinite quadratic with linear
costs optimal control of Markov jump with multiplicative noise systems,"
Automatica, vol. 43, pp. 587-597, 2007.

\bibitem{14} O. L. V. Costa and L. P. Wanderlei, ``Generalized Coupled Algebraic Riccati Equations for Discrete-time
Markov Jump with Multiplicative Noise Systems,"
European Journal of Control, vol. 5, pp. 391-408, 2008.

\bibitem{15}  M. A. Rami, X. Chen, J. B. Moore, and X. Zhou, ``Solvability
and asymptotic behavior of generalized Riccati equations
arising in indefinite stochastic LQ controls," IEEE Trans.
Automat. Control., vol.46, no.3, pp. 428-440, 2001.

\bibitem{16}   M. A. Rami and X. Zhou, ''Linear matrix inequalities, Riccati
equations, and indefinite stochastic linear quadratic control,"
IEEE Trans. Autom. Control., vol. 45, no. 6, pp. 1131-1142,
2000.

\bibitem{17}    O. L. V. Costa and M. D. Fragoso, ``Stability results for discrete-time
linear systems with Markovian jumping parameters," J. Math. Anal.
Appl., vol. 179, pp. 154-178, 1993.

\bibitem{18} H. Zhang, L. Li, J. Xu, and M. Fu, ``Linear quadratic regulation and stabilization of discrete-
time systems with delay and multiplicative noise," IEEE Trans. Autom. Control, vol. 60,
no. 10, pp. 2599-2613, 2015.

\bibitem{19} H. Zhang and J. Xu, ``Control for It$\hat{o}$ stochastic systems with
input delay," IEEE Trans. Autom. Control., vol. 62, no.1,
pp. 350-365, 2017.

\bibitem{20} Q. Qi, and H. Zhang, ``A Complete Solution to Optimal Control and Stabilization for Mean-field Systems: Part I, Discrete-time Case," arXiv preprint arXiv:
1608.06363, pp. 1-20, 2016.

\bibitem{21} M. Aitrami, X. Chen, and X. Zhou, ``Discrete-time indefinite LQ control
with state and control dependent noises," J. Global Optimizat., vol.
23, no. 3-4, pp. 245-265, 2002.

\bibitem{22}  A. Albert. ``Conditions for positive and nonnegative
definiteness in terms of pseudo-inverse," SIAM J. Appl. Math.,
vol.17, pp. 434-440, 1969.

\bibitem{23}  W. Zhang and B. S. Chen, ``On stabilizability and exact
observability of stochastic systems with their applications,"
Automatica, vol. 40, pp. 87-94, 2004.

\bibitem{24} H. Zhang, H. Wang, and L. Li, ``Adapted and casual maximum principle
and analytical solution to optimal control for stochastic multiplicative-
noise systems with multiple input-delays," in Proc. 51st IEEE Conf.
Decision Control, Maui, HI, USA,  pp. 2122-2127, 2012.

\end{thebibliography}

%

\end{document}